# INFORMATION BOUNDS FOR COX REGRESSION MODELS WITH MISSING DATA

By Bin Nan, Mary J. Emond[1] and Jon A. Wellner[2]

*University of Michigan and University of Washington*

We derive information bounds for the regression parameters in Cox models when data are missing at random. These calculations are of interest for understanding the behavior of efficient estimation in case-cohort designs, a type of two-phase design often used in cohort studies. The derivations make use of key lemmas appearing in Robins, Rotnitzky and Zhao [*J. Amer. Statist. Assoc.* **89** (1994) 846–866] and Robins, Hsieh and Newey [*J. Roy. Statist. Soc. Ser. B* **57** (1995) 409–424], but in a form suited for our purposes here. We begin by summarizing the results of Robins, Rotnitzky and Zhao in a form that leads directly to the projection method which will be of use for our model of interest. We then proceed to derive new information bounds for the regression parameters of the Cox model with data Missing At Random (MAR). In the final section we exemplify our calculations with several models of interest in cohort studies, including an i.i.d. version of the classical case-cohort design of Prentice [*Biometrika* **73** (1986) 1–11] and Self and Prentice [*Ann. Statist.* **16** (1988) 64–81].

**1. Introduction.** Models for missing data have been the subject of intense research over the past decade. In particular, the landmark paper of Robins, Rotnitzky and Zhao (1994) (hereafter RRZ) provides theoretical results for information bounds in semiparametric regression models with some covariates missing at random. RRZ studied extensively the special case where the model for the complete data is restricted only by specification of its mean, conditional on the covariates. They provided a brief treatment of the case where the full data model is the Cox regression model. In related

Received June 2002; revised January 2003.
[1]Supported in part by NIH Grant R29CA77607.
[2]Supported in part by NSF Grant DMS-9532039 and NIAID Grant 2R01 AI291968-04.

work, Robins, Hsieh and Newey (1995) (hereafter RHN) provided information bounds for classical regression models with missing covariate data.

Meanwhile, case-cohort and stratified case-cohort designs have become increasingly important and popular in epidemiology since the basic work of Prentice (1986) and Self and Prentice (1988). For reports of studies using these designs, see, for example, Bell, Hertz-Picciotto and Beaumont (2001), Dome, Chung, Bergemann, Umbricht, Saji, Carey, Grundy, Perlman, Breslow and Sukumar (1999), Margolis, Knauss and Bilker (2002), Mark, Qiao, Dawsey, Wu, Katki, Guntere, Fraumeni, Blot (2000), Rasmussen, Folsom, Catellier, Tsai, Garg and Eckfeldt (2001), Zeegers,Goldbohm and van den (2001) and Zeegers, Swaen, Kant, Goldbohm and van den Brandt (2001). These study designs correspond to missing data models, since complete data are collected only on a subsample of the study cohort. The currently used estimators for the Cox model with these designs are not known to be efficient, being based on pseudo-likelihoods or various ad hoc estimating equations. Because of the sheer volume of studies using these designs, it is becoming increasingly important to better understand the following:

(1) What are the information bounds for these types of designs and models?
(2) How much information is being lost by use of ad hoc estimators?
(3) Is it possible to construct reasonable, easily computable estimators which achieve the information bounds?

Our goal here is to begin to address the first two of these issues.

We begin by reorganizing and summarizing some results appearing in RRZ and RHN. Our summary (in Section 2) is formulated in a way which will lead quickly to information bounds for the models of primary concern here, namely Cox regression models with missing data. Our new information bounds for the Cox regression model with missing data are presented in Section 3. The efficient scores are characterized in terms of the solution of an integral equation.

In Section 4 the information bounds are calculated explicitly for particular submodels in several special cases, including case-cohort and exposure-stratified case-cohort versions of the Cox model. Although it has been known for some time that pseudo-likelihood estimators are not semiparametrically efficient, our explicit calculations quantify the loss of efficiency, and also show that two-phase designs with stratified subsampling can partially recover the information that is lost due to missing data.

Although we will not address question (3) in this paper, we note that for complex models such as those under study in Sections 3 and 4 of this paper, it is not uncommon for the calculation of information bounds to precede and aid in the development of efficient estimators. For example, the information bounds obtained by Sasieni (1992a, b) came seven or eight years before the development of efficient estimators for "partly linear" extensions of the Cox model in Huang (1999). Construction of efficient estimators for case-cohort designs, with and without stratification, will be treated by the first author. For preliminary work in this direction, see Nan (2001).

While our focus here is on information bounds rather than on construction of estimators, we comment briefly here on work on the estimation side of the problem. Most of the recent work on estimators for missing data in the Cox model focuses on improvements of the pseudo-likelihood estimators of Self and Prentice (1988); see, for example, Borgan, Langholz, Samuelsen, Goldstein and Pogoda (2000), Chen and Lo (1999) and the methods developed for related missing data models in Chatterjee, Chen and Breslow (2003).

**2. Information bounds for models with missing data.** We first give a brief review of the general setting for information bound calculations with missing data. This material is a reworking of important results in Robins, Rotnitzky and Zhao (1994) in a form suitable for our present calculations. Readers new to these calculations may also be interested in van der Vaart [(1998), pages 379–383] and Emond and Wellner (1995).

The general setup in this article is as follows: we suppose that $U^0$ is a random vector with distribution $Q$ in the model $\mathcal{Q}$: $U^0$ represents the "full" or "complete data." The complete or "full data" model $\mathcal{Q}$ may be parametric, semiparametric or nonparametric, but in our examples it will be semiparametric: $\mathcal{Q} = \{Q_{\theta,\eta} : \theta \in \Theta \subset \mathbb{R}^d, \eta \in \mathcal{H}\}$ where $\theta$ is the parameter of primary interest and $\eta$ is an infinite-dimensional "nuisance parameter." The "observed data" is $U$, where typically $U^0 = (U_1^0, U_2^0)$, and then $U = (U^0, R) = (U^0, 1)$ when the indicator variable $R = 1$, and $U = (U_1^0, R) = (U_1^0, 0)$ when $R = 0$. The distribution of $U$ is $P$, an element of the (induced) "observed data" model $\mathcal{P}$. In our examples $\mathcal{P}$ is semiparametric, parametrized by $(\theta, \eta)$, where $\theta \in \Theta \subset \mathbb{R}^d$ is the parameter of interest and $\eta$ is a nuisance parameter. The goal is to find the information bound for estimation of $\theta$ when $\eta$ is unknown based on observation of $U_1, \ldots, U_n$ i.i.d. as $U \sim P_{\theta,\eta} \in \mathcal{P}$.

Here is the primary model of interest for which the information bound is derived in Sections 3–5.

EXAMPLE (The Cox model with missing covariates). Let $T$ be a failure time, $C$ be a censoring time and $Z = (X, V) \in \mathbb{R}^d$ be a covariate vector which is not time dependent. The data $X$ are missing at random, while $Y \equiv T \wedge C$, $\Delta \equiv \mathbb{1}_{[T \leq C]}$ and $V$ and $R$ are always observed, where $R$ is an indicator of missingness as above. The full data are $U^0 = (Y, \Delta, X, V)$ and the observed data are $U = (Y, \Delta, RX, V, R)$ in the general notation introduced above. Note that $X$ may be missing by design, as in two-phase studies. In a two-phase study $(Y, \Delta, V)$ is observed for all subjects in phase 1 of the study. [In some "classical" case-cohort designs, only $(\Delta, V)$ is observed in phase 1. We will treat this case briefly in Section 3.3.] In phase 2 $X$ is obtained on a subsample of the subjects. The probability of being included in this subsample may depend on what was observed in phase 1. We are interested in estimating the effect on $T$ of the covariate $Z = (X, V)$. Let $(T|Z) \sim F(\cdot|Z)$ with density $f = f_{\theta,\lambda}$, where

$$1 - F_{\theta,\lambda}(t|z) = \exp(-e^{\theta'z}\Lambda(t)) \quad \text{so} \quad \frac{f_{\theta,\lambda}(t|z)}{1 - F_{\theta,\lambda}(t|z)} = e^{\theta'z}\lambda(t),$$

where $\lambda$ is the (Lebesgue) density of $\Lambda$. Also, let $(C|Z) \sim G(\cdot|Z)$ with density $g$, where

$$\frac{g(c|z)}{1 - G(c|z)} = \lambda_G(c|z) \quad \text{and} \quad \Lambda_G(c|z) = \int_0^c \lambda_G(t|z)\,dt.$$

We assume that $T$ and $C$ are conditionally independent given $Z$ (noninformative censoring). Let $Z \sim H$ with density $h$. Then $\mathcal{Q}$ is the set of all densities of the form

$$
\begin{aligned}
q_{\theta,\lambda,\lambda_G,h}&(y,\delta,z) \\
&= q(y,\delta,x,v) \\
(2.1) \quad &= \left[f(y|z)\int_{(y,\infty)} g(t|z)\,dt\right]^{\delta}\left[g(y|z)\int_{(y,\infty)} f(t|z)\,dt\right]^{1-\delta} h(z) \\
&= (e^{\theta'z}\lambda(y))^{\delta}\exp(-e^{\theta'z}\Lambda(y))(\lambda_G(y|z))^{1-\delta}\exp(-\Lambda_G(y|z))h(z).
\end{aligned}
$$

The regression parameter $\theta$ is the parameter of interest, and the nuisance parameter $\eta$ is $(\lambda, \lambda_G, h)$. This is basically the model introduced by Cox (1972). The model $\mathcal{P}$ for the observed data is the set of all distributions with densities of the form

$$
\begin{aligned}
p(r,y,\delta,(r\cdot x),v) &= (\pi(y,\delta,v)q(y,\delta,x,v))^r \\
(2.2) \quad &\times \left((1-\pi(y,\delta,v))\int q(y,\delta,x,v)\,d\mu(x)\right)^{1-r},
\end{aligned}
$$

where $\pi(Y,\Delta,V) = P(R=1|U_1^0)$, with $\pi(Y,\Delta,V) \geq \sigma > 0$, and $\mu$ is a dominating measure on $\mathcal{X}$.

Now we give a brief introduction to efficiency calculations in general; for more details see Bickel, Klaassen, Ritov and Wellner (1993) or van der Vaart [(1998), pages 362–371].

2.1. *Introduction to information bounds for semiparametric models.* The information bound for estimation of $\theta$ in the model $\mathcal{P}$ is equal to $I_\theta^{*-1}$. Here, $I_\theta^*$ is the efficient information matrix for $\theta$ in $\mathcal{P}$, given by

$$(2.3) \qquad I_\theta^* = E_P(l_\theta^* l_\theta^{*T}) \equiv E_P(l_\theta^{*\otimes 2}),$$

where $l_\theta^*$ is the efficient score for $\theta$. The efficient score for $\theta$ in a model $\mathcal{P} = \{P_{\theta,\eta} : \theta \in \Theta, \eta \in \mathcal{H}\}$ with nuisance parameters $\eta$ is the (componentwise) projection of the vector of scores $\dot{l}_\theta \in (L_2^0(P))^d$ for $\theta$ onto the orthogonal complement of the (closure of the linear span) of all scores for the nuisance parameters, $\dot{\mathcal{P}}_\eta$. Intuitively, when $\eta$ is unknown information about $\theta$ can only come from that component of $\dot{l}_\theta$ that is statistically independent of variability in the data controlled by the nuisance parameter. This component is $l_\theta^*$. Formally, each component of the efficient score for $\theta$ is orthogonal to all scores for nuisance parameters, where orthogonality is relative to the inner product $\langle b(U), a(U) \rangle_{L_2^0(P)} \equiv E_P\{ba\}$ in the space $L_2^0(P)$, the space of all mean-zero square-integrable functions of $U$. Let $\dot{\mathcal{P}}$ be the tangent space for $\mathcal{P}$. $\dot{\mathcal{P}}$ is the closure of the linear span of the scores of all submodels of $\mathcal{P}$ passing through $P$ (see BKRW). The tangent space for a model can be thought of as the space of all components of its score, and the "size" of $\dot{\mathcal{P}}$ corresponds to the amount of unknown information about $\mathcal{P}$. For example, when $\mathcal{P}$ is completely nonparametric, $\dot{\mathcal{P}}$ is all of $L_2^0(P)$. Now let $\dot{\mathcal{P}}_\theta$ and $\dot{\mathcal{P}}_\eta$ be the tangent spaces for the submodels of $\mathcal{P}$ where $\eta$ and $\theta$ are assumed known, respectively. Analogous definitions hold for $\dot{\mathcal{Q}}_\theta$ and $\dot{\mathcal{Q}}_\eta$. Then $\dot{\mathcal{P}}_\theta + \dot{\mathcal{P}}_\eta \subset \dot{\mathcal{P}}$ and we may assume for our purposes that $\dot{\mathcal{P}}_\theta + \dot{\mathcal{P}}_\eta = \dot{\mathcal{P}}$; see BKRW (page 76) for a discussion. The orthogonality condition described above for $l_\theta^* \in (\dot{\mathcal{P}}_\eta^\perp)^d$ is

$$(2.4) \qquad l_\theta^* \perp \dot{\mathcal{P}}_\eta \qquad \text{in } L_2^0(P);$$

that is, $E_P(l_\theta^* b) = 0$ for all $b \in \dot{\mathcal{P}}_\eta$, where the orgonality is componentwise (i.e., it holds for each component of the vector of functions $l_\theta^*$). Thus, our approach to calculating $l_\theta^*$ and the resulting efficiency bound is to project $\dot{l}_\theta$ onto the orthocomplement $\dot{\mathcal{P}}_\eta^\perp$ of $\dot{\mathcal{P}}_\eta$ in $L_2^0(P)$:

$$(2.5) \qquad l_\theta^* = \mathbf{\Pi}(\dot{l}_\theta | \dot{\mathcal{P}}_\eta^\perp),$$

where $\mathbf{\Pi}$ denotes the projection operator; see, for example, Bickel, Klaassen, Ritov and Wellner [(1993), Appendix A.2]. [Here the $\perp$ (orthogonal complement) denotes the ortho-complement in $L_2^0(P)$ or $L_2^0(Q)$, depending on the context.] The space $\dot{\mathcal{P}}_\eta^\perp$ is of further importance, because it contains all influence functions for all regular estimators of $\theta$ in $\mathcal{P}$. Note that $l_\theta^* = \mathbf{\Pi}(\dot{l}_\theta | \dot{\mathcal{P}}_\eta^\perp) = b^*$ if and only if

$$(2.6) \qquad \langle b, \dot{l}_\theta - b^* \rangle = 0 \qquad \text{for all } b \in \dot{\mathcal{P}}_\eta^\perp.$$

This implies that we can find the desired projection by proposing a guess $b^*$ for $l_\theta^*$ and then showing that (2.6) holds. This requires some understanding of $\dot{\mathcal{P}}_\eta^\perp$. However, this last requirement can be relaxed somewhat. Since $\dot{l}_\theta \in \dot{\mathcal{P}} = \dot{\mathcal{P}}_\theta + \dot{\mathcal{P}}_\eta$, we have

$$\mathbf{\Pi}(\dot{l}_\theta | \dot{\mathcal{P}}_\eta^\perp) = \mathbf{\Pi}(\dot{l}_\theta | \dot{\mathcal{P}} \cap \dot{\mathcal{P}}_\eta^\perp) = \mathbf{\Pi}(\dot{l}_\theta | \mathcal{M} \cap \dot{\mathcal{P}}_\eta^\perp)$$

for any (closed) subspace $\mathcal{M}$ such that $\dot{\mathcal{P}} \subset \mathcal{M} \subset L_2^0(P)$. So it is sufficient to be able to identify all $b$ in some set $\mathcal{M} \cap \dot{\mathcal{P}}_\eta^\perp$, which might not be all of $\dot{\mathcal{P}}_\eta^\perp$. Note that $\dot{\mathcal{P}} \cap \dot{\mathcal{P}}_\eta^\perp \subset \mathcal{M} \cap \dot{\mathcal{P}}_\eta^\perp \subset \dot{\mathcal{P}}_\eta^\perp$. This is essentially the approach of Robins, Rotnitzky and Zhao (1994); see also the discussion in van der Vaart [(1998), pages 379–383]. The approach proceeds by identifying an "intermediate" set $\mathcal{K} = \mathcal{M} \cap \dot{\mathcal{P}}_\eta^\perp$, and then knowing $l_\theta^* \in (\mathcal{K})^d$ provides a general form for $l_\theta^*$. An expression for $l_\theta^*$ will be obtained by finding the specific element of $\mathcal{K}$ that is the projection of $\dot{l}_\theta$ onto $\mathcal{K}$. The next subsection provides the formal results necessary to carry out this approach.

2.2. *Information bounds with missing data.* The following material is essentially a special case of results appearing in Robins, Rotnitzky and Zhao (1994). Throughout this subsection we take the complete data to be $U^0 = (U_1^0, U_2^0) \sim Q \in \mathcal{Q}$, with $U_2^0$ Missing At Random (MAR) in the observed data: thus $R$ is an indicator of whether $U_2^0$ is observed with $P(R = 1|U^0) = P(R = 1|U_1^0) \equiv \pi(U_1^0)$ with $\pi(U_1^0) \geq \sigma > 0$. If $\pi = \pi_\gamma$ is allowed to be unknown via parameters $\gamma \in \Gamma$, then we let $\mathcal{R} = \{\pi_\gamma : \gamma \in \Gamma\}$ denote the model for the "missingness probabilities" $\pi$. The observed data are $U \equiv (U_1^0, U_2^0, R)^R (U_1^0, R)^{1-R} \equiv (U_1^0, RU_2^0, R)$. Because there is a measurable map from $(U^0, R)$ to $U$, the tangent space for $\mathcal{Q} \times \mathcal{R}$ can be mapped to $\dot{\mathcal{P}}$ via an operator $\mathbf{A}$ that we call the *score operator*. The tangent space $\dot{\mathcal{Q}} \times \dot{\mathcal{R}}$ for $\mathcal{Q} \times \mathcal{R}$ is just $\dot{\mathcal{Q}}$ under our assumption that $\pi(U_1^0)$ is known, and we indeed impose this assumption throughout the rest of this paper.

REMARK. We have not included $R$ in $U^0$ because the functions in $\dot{\mathcal{Q}}$ do not depend on $R$. However, if $\pi$ were partially unknown, then the parameters of $\pi$ would be additional nuisance parameters, $R$ would be included in $U^0$, and $\dot{\mathcal{Q}}$ would be replaced by $\dot{\mathcal{Q}} + \dot{\mathcal{R}}$ in Lemmas 2.1 and 2.2.

For $a \in L_2^0(Q)$ define the (score) operator $\mathbf{A} : L_2^0(Q) \to L_2^0(P)$ by

$$\mathbf{A}a(U) \equiv (\mathbf{A}a)(U) \equiv E\{a(U^0)|U\} = Ra(U^0) + (1-R)E(a(U^0)|U_1^0).$$

LEMMA 2.1.

A. $\{\mathbf{A}a(U) : a \in \dot{\mathcal{Q}}\} \subset \dot{\mathcal{P}}$.
B. *The adjoint* $\mathbf{A}^T : L_2^0(P) \to L_2^0(Q)$ *of* $\mathbf{A}$ *is given by* $\mathbf{A}^T b(U^0) \equiv E\{b(U)|U^0\}$ *for* $b \in L_2^0(P)$.
C. $\mathbf{A}^T \mathbf{A} a = \pi(U_1^0) a(U^0) + (1 - \pi(U_1^0)) E[a|U_1^0]$ *for* $a \in L_2^0(Q)$.
D. *The operator* $(\mathbf{A}^T \mathbf{A})^{-1}$ *is given by*

$$(\mathbf{A}^T \mathbf{A})^{-1} a = \frac{a(U^0)}{\pi(U_1^0)} - \frac{1 - \pi(U_1^0)}{\pi(U_1^0)} E[a|U_1^0].$$

Robins, Rotnitzky and Zhao (1994) denote the operator $\mathbf{A}^T\mathbf{A}$ by $\mathbf{m}$. Then C and D are special cases of their Propositions 8.2a2 and 8.2e. The next lemma is key since it identifies $\dot{\mathcal{P}}_\eta^\perp$ once $\dot{\mathcal{Q}}_\eta^\perp$ is known.

LEMMA 2.2. *Suppose that $\pi(U_1^0) \geq \sigma > 0$. Then:*

A. *Range($\mathbf{A}$) is closed [and so is Range($\mathbf{A}|_{\dot{\mathcal{Q}}}$) or the range of $\mathbf{A}$ restricted to any other closed subspace of $L_2^0(Q)$].*
B. *Let $b \in L_2^0(P)$. Then $b \in \dot{\mathcal{P}}_\eta^\perp$ if and only if $\mathbf{A}^T b \in \dot{\mathcal{Q}}_\eta^\perp$.*

Part B of Lemma 2.2 is Lemma A.6 of Robins, Rotnitzky and Zhao (1994), while part A of Lemma 2.2 is proved in Robins and Rotnitzky [(1992), proof of Theorem 4.1, page 326]. It is also in the proof of Lemma A.4 in Robins, Rotnitzky and Zhao [(1994), page 862]. Note that the condition of Lemma 2.2 is used by Robins, Rotnitzky and Zhao (1994) in the proofs of both parts A and B since RRZ's proof of Lemma A.6 uses their Lemma A.2 which in turn has their (36) as a hypothesis.

Now suppose that $\dot{\mathcal{Q}}_\eta^\perp$ is known. Then for the subspace $\mathcal{M}$ in the discussion earlier in this section we will take $\dot{\mathcal{P}}_\eta + \mathcal{K}$, where $\mathcal{K}$ consists of the closed subspace of all functions $k(U)$ of the form

$$k(U) = \frac{R}{\pi}\zeta(U^0) - \frac{R-\pi}{\pi}E[\zeta(U^0)|U_1^0],$$

where $\zeta(U^0) \in \dot{\mathcal{Q}}_\eta^\perp$. Moreover, let $\mathcal{J}$ be the set of all $j(U)$ such that

$$j(U) = \frac{R}{\pi}\zeta(U^0) - \frac{R-\pi}{\pi}\phi(U_1^0),$$

where $\zeta(U^0) \in \dot{\mathcal{Q}}_\eta^\perp$ and $\phi(U_1^0)$ is any function in $L_2(P_{U_1^0})$. The set $\mathcal{J}$ is discussed by Robins, Rotnitzky and Zhao (1994) and by van der Vaart [(1998), page 383]. As noted by RRZ, the particular function $\phi(U_1^0)$ in the definition of $k$ yields the smallest variance for a given function $\zeta$.

The next three lemmas and two propositions characterize $\mathcal{J}$ and $\mathcal{K}$ in terms of $\dot{\mathcal{P}}$ and $\dot{\mathcal{P}}_\eta^\perp$ and show that $\mathcal{K}$ has the desired properties. These propositions form the basis for our specific information bound calculations for the Cox model in the sections to follow.

The next lemma shows that every $b = \mathbf{A}a \in L_2^0(P)$ for $a \in L_2^0(Q)$ can be decomposed into the form $(R/\pi)(\mathbf{A}^T\mathbf{A})a - \mathbf{\Pi}((R/\pi)(\mathbf{A}^T\mathbf{A})a|\mathcal{J}^{(2)})$, where $\mathcal{J}^{(2)}$ is the subspace of $L_2^0(P)$ with form of the second term in the definition of the class $\mathcal{J}$. The following Lemma 2.3 is a special case of Lemma A.3 of Robins, Rotnitzky and Zhao (1994).

LEMMA 2.3. *Suppose $a(U^0) \in L_2^0(Q)$. Then $\mathbf{A}a \perp \frac{R-\pi}{\pi}\phi(U_1^0)$ for all $a \in L_2^0(Q)$, $\phi \in L_2(Q)$; equivalently*

$$\mathbf{A}a = \frac{R}{\pi}(\mathbf{A}^T\mathbf{A})(a) - \mathbf{\Pi}\Big(\frac{R}{\pi}(\mathbf{A}^T\mathbf{A})(a)\Big|\mathcal{J}^{(2)}\Big), \tag{2.7}$$

*where $\mathcal{J}^{(2)} \equiv \{(R/\pi - 1)\phi(U_1^0) : \phi(U_1^0) \in L_2(Q)\}$.*

PROPOSITION 2.1. $\mathcal{K} \subset \mathcal{J} \subset \dot{\mathcal{P}}_\eta^\perp$ *and* $\dot{\mathcal{P}}_\eta^\perp \cap \dot{\mathcal{P}} \subset \mathcal{K}$.

REMARK. Any function $b \in (\dot{\mathcal{P}}_\eta^\perp)^d$ with $\langle b, l_\theta^* \rangle = J$, the $d \times d$ identity matrix, is an influence function for estimation of $\theta \in \Theta \subset \mathbb{R}^d$ in the model $\mathcal{P}$ for the observed data (see, e.g., BKRW, page 73, Proposition 3.4.2), and the decomposition of $b = \mathbf{A}a$ given in (2.7) shows how $b = \mathbf{A}a \in (\dot{\mathcal{P}}_\eta^\perp)^d$ is related to the influence function $\mathbf{A}^T\mathbf{A}a$ for estimation of $\theta$ in the model $\mathcal{Q}$ for the complete data. Note that $(R/\pi)(\mathbf{A}^T\mathbf{A})a$ is then the influence function of an inverse-weighted estimator of $\theta$ in $\mathcal{P}$ of the basic type proposed by Horvitz and Thompson (1952).

LEMMA 2.4. *Given a subspace $\mathcal{M}$ of $L_2^0(P)$ there is a unique projection map $\mathbf{\Pi}(\cdot|\mathcal{M})$ from $L_2^0(P)$ onto $\mathcal{M}$. In particular:*

A. *For $a^* \in \dot{\mathcal{P}}_\eta$, $h \in L_2^0(P)$ we have $a^* = \mathbf{\Pi}(h|\dot{\mathcal{P}}_\eta)$ if and only if $\langle h - a^*, a \rangle = 0$ for all $a \in \dot{\mathcal{P}}_\eta$.*

B. *For $b^* \in \dot{\mathcal{P}}_\eta^\perp$, $h \in L_2^0(P)$ we have $b^* = \mathbf{\Pi}(h|\dot{\mathcal{P}}_\eta^\perp)$ if and only if $\langle b, h - b^* \rangle = 0$ for all $b \in \dot{\mathcal{P}}_\eta^\perp$.*

C. *For $b^* \in \dot{\mathcal{P}}_\eta^\perp \cap \dot{\mathcal{P}}$, $h \in \dot{\mathcal{P}}$, we have $b^* = \mathbf{\Pi}(h|\dot{\mathcal{P}}_\eta^\perp \cap \dot{\mathcal{P}})$ if and only if $\langle b, h - b^* \rangle = 0$ for all $b \in \dot{\mathcal{P}}_\eta^\perp \cap \dot{\mathcal{P}}$.*

D. *Suppose that $h \in \mathcal{M}$ with $\dot{\mathcal{P}} \subset \mathcal{M} \subset L_2^0(P)$, $\mathcal{M}$ a closed subspace. Then for $b^* \in \dot{\mathcal{P}}_\eta^\perp \cap \mathcal{M}$, $b^* = \mathbf{\Pi}(h|\dot{\mathcal{P}}_\eta^\perp \cap \mathcal{M})$ if and only if $\langle b, h - b^* \rangle = 0$ for all $b \in \dot{\mathcal{P}}_\eta^\perp \cap \mathcal{M}$.*

E. *For $h \in \dot{\mathcal{P}}$, the projection $\mathbf{\Pi}(h|\dot{\mathcal{P}}_\eta^\perp \cap \mathcal{M}) = \mathbf{\Pi}(h|\dot{\mathcal{P}}_\eta^\perp \cap \dot{\mathcal{P}}) \in \dot{\mathcal{P}}$.*

The following proposition is an immediate consequence of Proposition 2.1 and Lemma 2.4, part E. It is a special case of Proposition 8.1e1 of Robins, Rotnitzky and Zhao (1994).

PROPOSITION 2.2. $l_\theta^* = \mathbf{\Pi}[\dot{l}_\theta|\mathcal{K}] = \frac{R}{\pi}\zeta^* - \frac{R-\pi}{\pi}E[\zeta^*|U_1^0]$, *where $l_\theta^*$ is the efficient score for $\theta$ in model $\mathcal{P}$ and $\zeta^*$ is the unique element of $(\dot{\mathcal{Q}}_\eta^\perp)^d$ satisfying*

$$\mathbf{\Pi}\Big(\frac{1}{\pi}\zeta^* - \frac{1-\pi}{\pi}E[\zeta^*|U_1^0]\Big|\dot{\mathcal{Q}}_\eta^\perp\Big) = l_\theta^{*0}, \tag{2.8}$$

*where $l_\theta^{*0}$ is the efficient score for $\theta$ in the complete data model $\mathcal{Q}$.*

**3. The Cox model with missing covariate data.** In this section we discuss the efficient score and information calculations for Cox regression models with missing covariates as introduced in the example in Section 2. This model is very useful in epidemiology studies, especially in two-stage designs (also called two-phase designs) where the probabilities $\pi$ are determined by the investigator. Equation (2.1) gives the joint density of the complete data $(Y, \Delta, X, V) \equiv (Y, \Delta, Z)$ and (2.2) gives the joint density of the observed data $(Y, \Delta, RX, V, R)$. The finite-dimensional parameter $\theta$ is the parameter of interest and the nuisance parameter $\eta = (\lambda, \lambda_G, h)$ is a vector of three infinite-dimensional nuisance parameters. We will use Proposition 2.2 to obtain the efficient score $l_\theta^*$ in the observed model $\mathcal{P}$ for Cox regression. The efficient score function in the full model $\mathcal{Q}$, $l_\theta^{*0}$, has been studied by many authors, such as Efron (1977), Andersen and Gill (1982), Begun, Hall, Huang and Wellner (1983) and Bickel, Klaassen, Ritov and Wellner (1993). Hence our main job will be to characterize the space $\dot{\mathcal{Q}}_\eta^\perp$.

3.1. *The nuisance parameter tangent space $\dot{\mathcal{Q}}_\eta^\perp$ of the full data model $\mathcal{Q}$.* The scores for the parameter of interest $\theta$ and the score operators for the nuisance parameters $\lambda$, $\lambda_G$ and $h$ in the "full" model $\mathcal{Q}$ given in (2.1) are the following:

$$(3.1) \quad \dot{l}_1^0(Y, \Delta, Z) \equiv \dot{l}_\theta^0(Y, \Delta, Z) = \Delta Z - Z\Lambda(Y)e^{\theta' Z} = \int Z\, dM(t),$$

$$(3.2) \quad \begin{aligned} \dot{l}_2^0(Y, \Delta, Z) &\equiv \dot{l}_\lambda^0 a(Y, \Delta, Z) \\ &= \Delta a(Y) - e^{\theta' Z}\int_0^Y a(t)\, d\Lambda(t) = \int a(t)\, dM(t), \end{aligned}$$

$$(3.3) \quad \begin{aligned} \dot{l}_3^0(Y, \Delta, Z) &\equiv \dot{l}_{\lambda_G}^0 b(Y, \Delta, Z) \\ &= (1-\Delta)b(Y, Z) - \int_0^Y b(t, Z)\, d\Lambda_G(t|Z) \\ &= \int b(t, Z)\, dM_G(t), \end{aligned}$$

$$(3.4) \quad \dot{l}_4^0(Y, \Delta, Z) \equiv \dot{l}_h^0 c(Y, \Delta, Z) = c(Z),$$

where $M$ and $M_G$ are martingales, conditional on $Z$, for the failure and censoring counting processes, respectively;

$$M(t) = \Delta \mathbb{1}(Y \leq t) - \int_0^t \mathbb{1}(Y \geq s)\, d\Lambda(s|Z),$$

$$(3.5) \quad M_G(t) = (1-\Delta)\mathbb{1}(Y \leq t) - \int_0^t \mathbb{1}(Y \geq s)\, d\Lambda_G(s|Z),$$

$$a(t) = \frac{\partial}{\partial \chi} \log \lambda_\chi(t),$$

$$b(t, Z) = \frac{\partial}{\partial \psi} \log \lambda_{G\psi}(t|Z), \qquad c(Z) = \frac{\partial}{\partial \kappa} \log h_\kappa(Z),$$

for regular parametric submodels $\{\lambda_\chi\}$, $\{\lambda_{G\psi}\}$ and $\{h_\kappa\}$ passing through the true parameters $\lambda$, $\lambda_G$ and $h$ when $\chi = 0$, $\psi = 0$ and $\kappa = 0$, respectively. Here we abuse notation slightly by writing $\Lambda(\cdot)$ for the baseline cumulative hazard function and $\Lambda(\cdot|Z)$ for the cumulative hazard function conditional on $Z$. Under the proportional hazards assumption $\Lambda(\cdot|Z) = \Lambda(\cdot)e^{\theta' Z}$.

Then the scores in the observed model (2.2) are the following:

$$\dot{l}_i(Y, \Delta, RX, V, R) = R\dot{l}_i^0(Y, \Delta, Z) + (1-R)E\{\dot{l}_i^0(Y, \Delta, Z)|Y, \Delta, V\},$$
$$i = 1, \ldots, 4.$$

Hence the tangent spaces for the parameters in the two models are

(3.6) $$\dot{\mathcal{Q}}_i \equiv [\dot{l}_i^0], \qquad \dot{\mathcal{P}}_i \equiv [\dot{l}_i], \qquad i = 1, \ldots, 4.$$

Here $[\alpha]$ denotes the closed linear span of the set $\alpha$ in $L_2^0(Q)$ or $L_2^0(P)$ as in Bickel, Klaassen, Ritov and Wellner [(1993), page 49]. By definition all the elements in $\dot{\mathcal{Q}}_i$ and $\dot{\mathcal{P}}_i$, $i = 1, \ldots, 4$, are square integrable. It is easy to see that $\dot{\mathcal{Q}}_2$, $\dot{\mathcal{Q}}_3$ and $\dot{\mathcal{Q}}_4$ are mutually orthogonal. Since they are closed (by definition), the nuisance tangent space is $\dot{\mathcal{Q}}_\eta = \dot{\mathcal{Q}}_2 + \dot{\mathcal{Q}}_3 + \dot{\mathcal{Q}}_4$.

Let $W_1$ and $W_2$ be subdistributions on $\mathbb{R}^{d+1}$ defined by

$$W_1(y, z) \equiv Q(Y \leq y, Z \leq z, \Delta = 1)$$
$$= \int_{(-\infty, z]} \int_{(0, y]} (1 - G(t|z'))\, dF(t|z')\, dH(z'),$$

$$W_2(y, z) \equiv Q(Y \leq y, Z \leq z, \Delta = 0)$$
$$= \int_{(-\infty, z]} \int_{(0, y]} (1 - F(t|z'))\, dG(t|z')\, dH(z'),$$

corresponding to the uncensored and censored data. Hence $W = W_1 + W_2$ is the marginal distribution of $(Y, Z)$. Then we can define $L_2$ spaces corresponding to the subprobability measures

$$L_2(W_i) = \left\{u(Y, Z) : \int\int u^2(y, z)\, dW_i(y, z) < \infty\right\}, \qquad i = 1, 2.$$

These spaces will be used in characterizing $\dot{\mathcal{Q}}_2$, $\dot{\mathcal{Q}}_3$ and thus $\dot{\mathcal{Q}}_\eta^\perp$. It is easy to see that $L_2(W) \equiv L_2(Q_{Y,Z}) = L_2(W_1) \cap L_2(W_2)$, $L_2(Q_{T,Z}) \subset L_2(W_1)$ and $L_2(Q_{C,Z}) \subset L_2(W_2)$. Here $Q_{Y,Z}$, $Q_{T,Z}$ and $Q_{C,Z}$ denote the marginal distributions of $(Y, Z)$, $(T, Z)$ and $(C, Z)$ under $Q \in \mathcal{Q}$.

Then the conditional distribution and conditional subdistributions of $Y$ given $Z$, $W(y|z)$, $W_1(y|z)$ and $W_2(y|z)$ have the following forms:

$$W(y|z) = \int_0^y (1 - G(t|z))\, dF(t|z) + \int_0^y (1 - F(t|z))\, dG(t|z)$$
$$= 1 - (1 - F(y|z))(1 - G(y|z)),$$
$$W_1(y|z) = \int_0^y (1 - G(t|z))\, dF(t|z),$$
$$W_2(y|z) = \int_0^y (1 - F(t|z))\, dG(t|z).$$

Now we define two operators $\mathbf{R}_1$ and $\mathbf{R}_2$ as follows:

$$(3.7) \quad \mathbf{R}_1 b(y, z) = b(y, 1, z) - \frac{\int_y^\infty b(t, 1, z)\, dW_1(t|z) + b(t, 0, z)\, dW_2(t|z)}{1 - W(y|z)},$$

$$(3.8) \quad \mathbf{R}_2 b(y, z) = b(y, 0, z) - \frac{\int_y^\infty b(t, 1, z)\, dW_1(t|z) + b(t, 0, z)\, dW_2(t|z)}{1 - W(y|z)}.$$

They can be rewritten as

$$(3.9) \qquad \mathbf{R}_1 b(y, z) = b(y, 1, z) - E[b(Y, \Delta, Z)|Y > y, Z = z],$$
$$(3.10) \qquad \mathbf{R}_2 b(y, z) = b(y, 0, z) - E[b(Y, \Delta, Z)|Y > y, Z = z].$$

We will show later in Proposition 3.1 that $\mathbf{R}_1$ and $\mathbf{R}_2$ map $L_2^0(Q)$ to $L_2(W_1)$ and $L_2(W_2)$, respectively. These operators are similar to the $\mathbf{R}$ operator discussed by Ritov and Wellner (1988), Efron and Johnstone (1990) and Bickel, Klaassen, Ritov and Wellner (1993).

The following Proposition 3.1 plays a key role in characterizing the space $\dot{\mathcal{Q}}_\eta^\perp$ which will be used to derive the efficient score $l_\theta^*$ for $\theta$ in the model $\mathcal{P}$ in the next subsection.

PROPOSITION 3.1. *Any function $b \in L_2^0(Q)$ can be decomposed as follows:*

$$(3.11) \quad b(Y, \Delta, Z) = \int \mathbf{R}_1 b(t, Z)\, dM(t) + \int \mathbf{R}_2 b(t, Z)\, dM_G(t) + E(b|Z),$$

*where*

$$(3.12) \qquad \mathbf{R}_1 b(Y, Z) \in L_2(W_1), \qquad \mathbf{R}_2 b(Y, Z) \in L_2(W_2).$$

*The decomposition is unique in the sense that $\mathbf{R}_1 b$ is unique a.e. $W_1$ and $\mathbf{R}_2 b$ is unique a.e. $W_2$.*

To prove Proposition 3.1, we will use the following two lemmas.

LEMMA 3.1. *For the failure counting process martingale $\{M(t): t \geq 0\}$ defined by (3.5), we have*

$$(3.13) \quad \int h_1(t, Z) \, dM(t) \in L_2^0(Q) \qquad \text{if and only if } h_1(Y, Z) \in L_2(W_1),$$

*and similarly for the censoring counting process martingale,*

$$(3.14) \quad \int h_2(t, Z) \, dM_G(t) \in L_2^0(Q) \qquad \text{if and only if } h_2(Y, Z) \in L_2(W_2).$$

PROOF. By using the methods of Chapter 6 of Shorack and Wellner (1986) we can easily show that

$$E\left(\int h_1(t, Z) \, dM(t)\right)^2 = \iint h_1^2(t, s) \, dW_1(t, s),$$

$$E\left(\int h_2(t, Z) \, dM_G(t)\right)^2 = \iint h_2^2(t, s) \, dW_2(t, s).$$

The zero means are trivial. □

LEMMA 3.2. *For any functions $h_j : \mathbb{R}^+ \times \mathbb{R}^d \to \mathbb{R}^q$ in $L_2(W_1)$, $j = 1, 2$,*

$$E\left[\int h_1(t, Z) \, dM(t) \int h_2(t, Z) \, dM(t)\right] = E[\Delta h_1(Y, Z) h_2(Y, Z)].$$

*Similarly, for any functions $h_j : \mathbb{R}^+ \times \mathbb{R}^d \to \mathbb{R}^q$ in $L_2(W_2)$, $j = 1, 2$,*

$$E\left[\int h_1(t, Z) \, dM_G(t) \int h_2(t, Z) \, dM_G(t)\right] = E[(1 - \Delta) h_1(Y, Z) h_2(Y, Z)].$$

PROOF. This follows from Lemma 1 of Sasieni (1992a). □

PROOF OF PROPOSITION 3.1. Equation (3.11) can be verified directly by the definitions of $\mathbf{R}_1$ and $\mathbf{R}_2$ operators in (3.7) and (3.8). See Nan (2001) for details. By Lemma 3.1 we know that the right-hand side of (3.11) is in $L_2^0(Q)$ if (3.12) is true. Now we show (3.12). Let

$$m(y, z) = \frac{\int_y^\infty b(t, 1, z) \, dW_1(t|z) + b(t, 0, z) \, dW_2(t|z)}{1 - W(y|z)}.$$

Obviously, $b(Y, 1, Z) \in L_2(W_1)$ and $b(Y, 0, Z) \in L_2(W_2)$ since

$$E[b^2(Y, \Delta, Z)] = E[\Delta b^2(Y, 1, Z) + (1 - \Delta) b^2(Y, 0, Z)]$$
$$= \iint b^2(y, 1, z) \, dW_1(y, z) + \iint b^2(y, 0, z) \, dW_2(y, z) < \infty.$$

Thus we only need to show that $m(Y, Z) \in L_2(W)$. We rewrite $m(Y, Z)$ as

$$m(Y, Z) = \frac{1}{1 - W(Y|Z)} \int_Y^\infty \{b(t, 1, Z)\alpha(t, Z) + b(t, 0, Z)\beta(t, Z)\} \, dW(t|Z),$$

where

$$\alpha(t, Z) = \frac{dW_1(t|Z)}{dW(t|Z)}, \qquad \beta(t, Z) = \frac{dW_2(t|Z)}{dW(t|Z)}.$$

It is obvious that $0 \leq \alpha \leq 1$ and $0 \leq \beta \leq 1$. Then by the same argument as that on page 423 of Bickel, Klaassen, Ritov and Wellner (1993), we have

$E[m^2(Y, Z)|Z]$

$$\leq 4 \int \left\{ b(t, 1, Z) \frac{dW_1(t|Z)}{dW(t|Z)} + b(t, 0, Z) \frac{dW_2(t|Z)}{dW(t|Z)} \right\}^2 dW(t|Z)$$

$$\leq 8 \left\{ \int b^2(t, 1, Z)\alpha(t, Z) \, dW_1(t|Z) + \int b^2(t, 0, Z)\beta(t, Z) \, dW_2(t|Z) \right\}$$

$$\leq 8 \left\{ \int b^2(t, 1, Z) \, dW_1(t|Z) + \int b^2(t, 0, Z) \, dW_2(t|Z) \right\}.$$

Then by Fubini's theorem we have

$$E[m^2(Y, Z)] \leq 8 \left\{ \iint b^2(t, 1, s) \, dW_1(t, s) + \iint b^2(t, 0, s) \, dW_2(t, s) \right\} < \infty.$$

Actually, from the above proof we have also shown that

$$L_2^0(Q) \equiv \left\{ h_0(Z) + \int h_1(t, Z) \, dM(t) \right.$$

$$\left. + \int h_2(t, Z) \, dM_G(t) : h_0 \in L_2^0(H), h_1 \in L_2(W_1), h_2 \in L_2(W_2) \right\}.$$

The uniqueness can be proved by showing that any two decompositions of an element in $L_2^0(Q)$ are identical. Suppose we have

$$b = \int h_1(t, Z) \, dM(t) + \int h_2(t, Z) \, dM_G(t) + h_0(Z);$$

$$b = \int h_1'(t, Z) \, dM(t) + \int h_2'(t, Z) \, dM_G(t) + h_0'(Z).$$

Taking the expectation of the square of the difference of the right-hand sides of the two equalities and using the orthogonality of $M$, $M_G$, and any function of $Z$, by Lemma 3.2 we have

$$0 = \int (h_1 - h_1')^2 \, dW_1 + \int (h_2 - h_2')^2 \, dW_2 + \int (h_0 - h_0')^2 \, dH.$$

Thus $h_0' = h_0$ a.s. $H$, $h_1' = h_1$ a.e. $W_1$ and $h_2' = h_2$ a.e. $W_2$. □

Now we are ready to discuss the space $\dot{\mathcal{Q}}_\eta^\perp$.

PROPOSITION 3.2. *For any function $s(Y,Z) \in L_2(W_1)$ define the operator $\mathbf{B}$ by*

$$(3.15) \qquad \mathbf{B}s = \int \{s(t,Z) - E[s(Y,Z)|Y=t, \Delta=1]\} \, dM(t).$$

*Then:*

(i) $\mathbf{B}s \perp \dot{\mathcal{Q}}_\eta$ *in* $L_2^0(Q)$.
(ii) *For any* $b \in L_2^0(Q)$ *we have* $\mathbf{\Pi}(b|\dot{\mathcal{Q}}_\eta^\perp) = \mathbf{B} \circ \mathbf{R}_1 b$.
(iii) $\dot{\mathcal{Q}}_\eta^\perp = (\dot{\mathcal{Q}}_2 + \dot{\mathcal{Q}}_3 + \dot{\mathcal{Q}}_4)^\perp = \{\mathbf{B}s : s \in L_2(W_1)\}$.

PROOF.  (i) Since $\dot{\mathcal{Q}}_2$, $\dot{\mathcal{Q}}_3$ and $\dot{\mathcal{Q}}_4$ are mutually orthogonal and $M \perp M_G$, we have

$$\mathbf{\Pi}(\mathbf{B}s|\dot{\mathcal{Q}}_2 + \dot{\mathcal{Q}}_3 + \dot{\mathcal{Q}}_4) = \mathbf{\Pi}(\mathbf{B}s|\dot{\mathcal{Q}}_2) + \mathbf{\Pi}(\mathbf{B}s|\dot{\mathcal{Q}}_3) + \mathbf{\Pi}(\mathbf{B}s|\dot{\mathcal{Q}}_4) = \mathbf{\Pi}(\mathbf{B}s|\dot{\mathcal{Q}}_2).$$

Let $m(Y,\Delta,Z) = \int s(t,Z) \, dM(t)$. Then $\mathbf{\Pi}(m|\dot{\mathcal{Q}}_2) = \int a^*(t) \, dM(t)$ for some $a^*(Y) \in L_2(W_1)$ satisfying

$$E\left[\int \{s(t,Z) - a^*(t)\} \, dM(t) \int a(t) \, dM(t)\right] = 0$$

for any $a(Y) \in L_2(W_1)$. Now by Lemma 3.2 the left-hand side of the above equation is equal to

$$E[\Delta\{s(Y,Z) - a^*(Y)\}a(Y)] = E[\{E[\Delta s(Y,Z)|Y] - a^*(Y)E[\Delta|Y]\}a(Y)]$$

and hence

$$a^*(Y) = \frac{E[\Delta s(Y,Z)|Y]}{E[\Delta|Y]} = E[s(Y,Z)|Y, \Delta=1].$$

So $\mathbf{B}s \perp \dot{\mathcal{Q}}_2$, and this yields $\mathbf{B}s \perp \dot{\mathcal{Q}}_2 + \dot{\mathcal{Q}}_3 + \dot{\mathcal{Q}}_4 \equiv \dot{\mathcal{Q}}_\eta$.

(ii) From Proposition 3.1 we know that for any $b \in L_2^0(Q)$ we have the decomposition (3.11). Hence, from the proof of part (i) we know that

$$\mathbf{\Pi}(b|\dot{\mathcal{Q}}_2 + \dot{\mathcal{Q}}_3 + \dot{\mathcal{Q}}_4)$$
$$= \mathbf{\Pi}(b|\dot{\mathcal{Q}}_2) + \mathbf{\Pi}(b|\dot{\mathcal{Q}}_3) + \mathbf{\Pi}(b|\dot{\mathcal{Q}}_4)$$
$$= \int E[\mathbf{R}_1 b(Y,Z)|Y=t, \Delta=1] \, dM(t)$$
$$\quad + \int \mathbf{R}_2 b(t,Z) \, dM_G(t) + E(b|Z).$$

Thus,

$$\mathbf{\Pi}(b|(\dot{\mathcal{Q}}_2 + \dot{\mathcal{Q}}_3 + \dot{\mathcal{Q}}_4)^\perp) = b - \mathbf{\Pi}(b|\dot{\mathcal{Q}}_2 + \dot{\mathcal{Q}}_3 + \dot{\mathcal{Q}}_4) = \mathbf{B} \circ \mathbf{R}_1 b.$$

(iii) This is an immediate consequence of (i) and (ii).  □

If we choose $s(t,Z) = Z$, then $\mathbf{B}s$ is the efficient score for $\theta$ in the ("full" or "complete" data) model $\mathcal{Q}$.

3.2. *Efficient score for $\theta$ in observed model $\mathcal{P}$.* In this subsection we will use the results in Section 2 and the previous subsection to derive the efficient score $l_\theta^*$ of $\mathcal{P}$. Since from Proposition 3.2 we know that $\dot{\mathcal{Q}}_\eta^\perp = \{\mathbf{B}s, s \in L_2(W_1)\}$, for the model $\mathcal{P}$, with $P \in \mathcal{P}$, we define the class $\mathcal{K}$ of all functions with the form

$$k(Y, \Delta, RX, V, R) = \frac{R}{\pi}\mathbf{B}(s(Y, X, V)) - \frac{R - \pi}{\pi}E[\mathbf{B}(s(Y, X, V))|Y, \Delta, V],$$

where $s \in L_2(W_1)$. Note that we can rewrite the functions $k$ in terms of the operator $\mathbf{D} \colon L_2(W_1) \to L_2^0(Q)$ defined by

(3.16) $$\mathbf{D}u(Y, \Delta, Z) = \int u(t, Z)\, dM(t).$$

Thus $\mathbf{B}s = \mathbf{D} \circ \mathbf{\Pi}_1 s = \mathbf{D}u$, where

(3.17) $$u(Y, Z) \equiv \mathbf{\Pi}_1 s \equiv s(Y, Z) - E[s(Y, Z)|Y, \Delta = 1].$$

Then we have the following proposition:

PROPOSITION 3.3. *The efficient score $l_\theta^*$ for $\theta$ in the model $\mathcal{P}$ for the observed data is given by*

$$l_\theta^*(R, Y, \Delta, R \cdot X, V)$$
$$= k^* \equiv \frac{R}{\pi}\mathbf{D}u^*(Y, \Delta, Z) - \frac{R - \pi}{\pi}E[\mathbf{D}u^*(Y, \Delta, Z)|Y, \Delta, V],$$

*where $u^* \in (L_2(W_1))^d$ is the unique a.e. $W_1$ solution of the equation*

(3.18) $$\begin{aligned}\mathbf{\Pi}_1 Z &= \mathbf{\Pi}_1 \circ \mathbf{R}_1\left\{\mathbf{D}u^* + \left(\frac{1-\pi}{\pi}\right)[\mathbf{D}u^* - E(\mathbf{D}u^*|Y, \Delta, V)]\right\} \\ &= u^* + \mathbf{T}u^*,\end{aligned}$$

*where*

(3.19) $$\mathbf{T} \equiv \mathbf{\Pi}_1 \circ \mathbf{R}_1 \circ \mathbf{H}$$

*and*

(3.20) $$\mathbf{H}u^* = \left(\frac{1-\pi}{\pi}\right)[\mathbf{D}u^* - E(\mathbf{D}u^*|Y, \Delta, V)].$$

COROLLARY 3.1. *The function $u^*$ in Proposition 3.3 also satisfies, equivalently,*

(3.21) $$\begin{aligned}u^*(Y, Z) &- \Big\{\mathbf{K}u^*(Y, Z) \\ &- \frac{\pi(Y, 1, V)}{E[\pi(Y, 1, V)|Y, \Delta = 1]}E[\mathbf{K}u^*(Y, Z)|TY, \Delta = 1]\Big\} \\ &= \pi(Y, 1, V)\{Z - E[Z|Y, R\Delta = 1]\},\end{aligned}$$

*where the operator* $\mathbf{K}$ *is a linear operator defined by*

$$
\begin{aligned}
\mathbf{K}u^*&(Y,Z)\\
&= -E[\mathbf{D}u^*(Y',\Delta,Z)|Y'>Y,Z]\\
&\quad + \pi(Y,1,V)E\left(\frac{\mathbf{D}u^*(Y',\Delta,Z)}{\pi(Y',\Delta,V)}\Big|Y'>Y,Z\right)\\
&\quad + (1-\pi(Y,1,V))E[\mathbf{D}u^*|Y,\Delta,V]|_{\Delta=1}\\
&\quad - \pi(Y,1,V)E\left(\frac{1-\pi(Y',\Delta,V)}{\pi(Y',\Delta,V)}\right.\\
&\qquad\qquad\qquad \times E[\mathbf{D}u^*(Y',\Delta,Z)|Y',\Delta,V]\Big|Y'>Y,Z\bigg).
\end{aligned}
\tag{3.22}
$$

Our proof of Proposition 3.3 will use the following lemma.

LEMMA 3.3. *Denote the adjoint of* $\mathbf{D}$ *by* $\mathbf{D}^T: L_2^0(Q) \to L_2(W_1)$. *Let* $\mathbf{\Pi}_1$ *be defined as in (*3.17*). Then:*

(i) $\mathbf{D}^T = \mathbf{R}_1$.
(ii) $\mathbf{R}_1 \circ \mathbf{D} = \mathbf{I}$.
(iii) $\mathbf{\Pi}_1$ *is a projection operator on* $L_2(W_1)$.

PROOF. (i) Let $a \in L_2(W_1)$, $b \in L_2^0(Q)$. Then we have $\langle a, \mathbf{D}^T b\rangle_{L_2(W_1)} = \langle \mathbf{D}a, b\rangle_{L_2^0(Q)}$. By Proposition 3.1,

$$b(Y,\Delta,Z) = \int \mathbf{R}_1 b(t,Z)\,dM(t) + \int \mathbf{R}_2 b(t,Z)\,dM_G(t) + E[b(Y,\Delta,Z)|Z].$$

Then we have

$$
\begin{aligned}
\langle \mathbf{D}a, b\rangle_{L_2^0(Q)} &= E_Q\bigg\{\int a(t,Z)\,dM(t)\int \mathbf{R}_1 b(t,Z)\,dM(t)\bigg\}\\
&= E_Q\{\Delta a(Y,Z)\mathbf{R}_1 b(Y,Z)\}\\
&= \iint a(t,z)\mathbf{R}_1 b(t,z)\,dW_1(t,z) = \langle a, \mathbf{R}_1 b\rangle_{L_2(W_1)}.
\end{aligned}
$$

(ii) For all $h \in L_2(W_1)$, by Proposition 3.1 we have $\mathbf{D}h = \int h\,dM = \int \mathbf{R}_1 \circ \mathbf{D}h\,dM$ and thus $h = \mathbf{R}_1 \circ \mathbf{D}h$ a.e. $W_1$.

(iii) That $\mathbf{\Pi}_2(\cdot) = E(\cdot|Y,\Delta=1)$ is a projection operator on $L_2(W_1)$ can be shown by checking the three properties in Proposition A.2.2 of Bickel, Klaassen, Ritov and Wellner (1993). So is $\mathbf{\Pi}_1 = \mathbf{I} - \mathbf{\Pi}_2$. $\square$

PROOF OF PROPOSITION 3.3.  We use Proposition 2.2 directly to prove Proposition 3.3. For the Cox regression model we have $l_\theta^{*0} = \mathbf{B}Z = \mathbf{D} \circ \mathbf{\Pi}_1 Z$ and $\zeta^* = \mathbf{D}u^*$. Let

$$b(Y, \Delta, Z) = \mathbf{B}Z - \frac{1}{\pi}\mathbf{D}u^* + \frac{1-\pi}{\pi}E[\mathbf{D}u^*|Y, \Delta, V].$$

Then by Proposition 2.2 and Proposition 3.2(ii) we have $\mathbf{\Pi}(b|\dot{\mathcal{Q}}_\eta^\perp) = \mathbf{B} \circ \mathbf{R}_1 b = 0$. Since by Lemma 3.3 we have $\mathbf{B} \circ \mathbf{R}_1 b = \mathbf{D} \circ \mathbf{\Pi}_1 \circ \mathbf{R}_1 b = \mathbf{D} \circ \mathbf{\Pi}_1 \mathbf{D}^T b = \mathbf{D} \circ (\mathbf{D} \circ \mathbf{\Pi}_1)^T b = \mathbf{D} \circ \mathbf{B}^T b$ and $\mathbf{R}_1 \circ \mathbf{D} = \mathbf{I}$, $\mathbf{B} \circ \mathbf{R}_1 b = 0$ implies $\mathbf{R}_1 \circ \mathbf{B} \circ \mathbf{R}_1 b = \mathbf{R}_1 \circ \mathbf{D} \circ \mathbf{B}^T b = \mathbf{B}^T b = 0$ which is

$$(3.23) \qquad \mathbf{B}^T\left(\mathbf{B}Z - \frac{1}{\pi}\mathbf{D}u^* + \frac{1-\pi}{\pi}E[\mathbf{D}u^*|Y, \Delta, V]\right) = 0.$$

Since $u^*(Y, Z) \equiv s^*(Y, Z) - E[s^*(Y, Z)|Y, \Delta = 1] \in L_2(W_1)$ satisfies

$$(3.24) \qquad E[u^*(Y, Z)|Y, \Delta = 1] = 0,$$

we must solve the pair of equations (3.23) and (3.24) for the function $u^*$. Note that $\mathbf{\Pi}_1$ is a projection operator and thus we have

$$\begin{aligned}\mathbf{B}^T \circ \mathbf{B}Z &= \mathbf{\Pi}_1^T \circ \mathbf{D}^T \circ \mathbf{B}Z \\ &= \mathbf{\Pi}_1 \circ \mathbf{R}_1 \circ \mathbf{D} \circ \mathbf{\Pi}_1 Z = \mathbf{\Pi}_1^2 Z = \mathbf{\Pi}_1 Z = Z - E(Z|Y, \Delta = 1).\end{aligned}$$

Hence (3.23) can be rewritten, using Lemma 3.3 (ii) and (iii), as

$$(3.25)\quad\begin{aligned}Z - E[Z|Y, \Delta = 1] &= \mathbf{\Pi}_1 \circ \mathbf{R}_1 \left\{\frac{\mathbf{D}u^*}{\pi} - \frac{1-\pi}{\pi}E(\mathbf{D}u^*|Y, \Delta, V)\right\} \\ &= \mathbf{\Pi}_1 \circ \mathbf{R}_1 \left\{\mathbf{D}u^* + \frac{1-\pi}{\pi}(\mathbf{D}u^* - E(\mathbf{D}u^*|Y, \Delta, V))\right\} \\ &= u^* + \mathbf{T}u^*,\end{aligned}$$

where $\mathbf{T} = \mathbf{\Pi}_1 \circ \mathbf{R}_1 \circ \mathbf{H}$ and $\mathbf{H}$ is given by (3.20). Thus (3.18) holds.

To see that the solution is unique, we argue as follows: let $\zeta^* = \mathbf{D}u^*$. Then from Proposition 2.2 we know that $\zeta^* \in (\dot{\mathcal{Q}}_\eta^\perp)^d$ is the unique solution of the operator equation

$$\mathbf{\Pi}\left(\frac{1}{\pi}\zeta^* - \frac{1-\pi}{\pi}E[\zeta^*|Y, \Delta, V]\Big|\dot{\mathcal{Q}}_\eta^\perp\right) = l_\theta^{*0}.$$

Suppose we have

$$\zeta^* = \mathbf{D}u_1^* = \int u_1^*(t, Z)\,dM(t) \quad \text{and} \quad \zeta^* = \mathbf{D}u_2^* = \int u_2^*(t, Z)\,dM(t).$$

Then taking the expectation of the square of the difference of the two equalities componentwise and using Lemma 3.2 (componentwise), we have
$$0 = \int |u_1^* - u_2^*|^2 \, dW_1.$$
It follows that $u_1^* = u_2^*$ a.e. $W_1$. □

PROOF OF COROLLARY 3.1. It also follows from (3.25) that

$$
\begin{aligned}
(3.26) \quad & Z - E[Z|Y, \Delta = 1] \\
&= \mathbf{\Pi}_1 \circ \mathbf{R}_1 \left\{ \frac{\mathbf{D}u^*}{\pi} - \frac{1-\pi}{\pi} E(\mathbf{D}u^*|Y, \Delta, V) \right\} \\
&= \mathbf{R}_1 \left\{ \frac{\mathbf{D}u^*}{\pi} - \frac{1-\pi}{\pi} E(\mathbf{D}u^*|Y, \Delta, V) \right\} \\
&\quad - E\left\{ \mathbf{R}_1 \left\{ \frac{\mathbf{D}u^*}{\pi} - \frac{1-\pi}{\pi} E(\mathbf{D}u^*|Y, \Delta, V) \right\} \Big| Y, \Delta = 1 \right\} \\
(3.27) \quad &\equiv \mathbf{R}_1 \left\{ \frac{\mathbf{D}u^*}{\pi} - \frac{1-\pi}{\pi} E(\mathbf{D}u^*|Y, \Delta, V) \right\} - f_{u^*}(Y).
\end{aligned}
$$

Now
$$
(3.28) \quad \mathbf{R}_1\left(\frac{\mathbf{D}u^*}{\pi}\right) = \frac{\mathbf{D}u^*|_{\Delta=1}}{\pi(Y, 1, V)} - E\left(\frac{\mathbf{D}u^*(Y', \Delta, Z)}{\pi(Y', \Delta, V)} \Big| Y' > Y, Z\right)
$$

and

$$
\begin{aligned}
& \mathbf{R}_1\left(\frac{1-\pi}{\pi} E[\mathbf{D}u^*|Y, \Delta, V]\right) \\
(3.29) \quad &= \frac{1 - \pi(Y, 1, V)}{\pi(Y, 1, V)} E[\mathbf{D}u^*|Y, \Delta, V]|_{\Delta=1} \\
&\quad - E\left(\frac{1 - \pi(Y', \Delta, Z)}{\pi(Y', \Delta, V)} E[\mathbf{D}u^*(Y', \Delta, Z)|Y', \Delta, V] \Big| Y' > Y, Z\right).
\end{aligned}
$$

Substituting (3.28) and (3.29) into (3.27) yields

$$
\begin{aligned}
& Z - E[Z|Y, \Delta = 1] + f_{u^*}(Y) \\
&= \frac{\mathbf{D}u^*|_{\Delta=1}}{\pi(Y, 1, V)} - E\left(\frac{\mathbf{D}u^*(Y', \Delta, Z)}{\pi(Y', \Delta, V)} \Big| Y' > Y, Z\right) \\
(3.30) \quad &\quad - \frac{1 - \pi(Y, 1, V)}{\pi(Y, 1, V)} E[\mathbf{D}u^*|Y, \Delta, V]|_{\Delta=1} \\
&\quad + E\left(\frac{1 - \pi(Y', \Delta, V)}{\pi(Y', \Delta, V)} E[\mathbf{D}u^*|Y', \Delta, V] \Big| Y' > Y, Z\right) \\
&\equiv \frac{1}{\pi(Y, 1, V)} (u^*(Y, Z) - \mathbf{K}u^*(Y, Z)).
\end{aligned}
$$

Here we use

$$\mathbf{D}u^*|_{\Delta=1} = \mathbf{R}_1 \circ \mathbf{D}u^*(Y,Z) + E[\mathbf{D}u^*(Y',\Delta,Z)|Y'>Y,Z],$$

and the operator $\mathbf{K}$ is defined as (3.22). Then we have

(3.31)  $u^*(Y,Z) - \mathbf{K}u^*(Y,Z) = \pi(Y,1,V)\{Z - E(Z|Y,\Delta=1) + f_{u^*}(Y)\}.$

Thus, taking conditional expectations given $Y, \Delta=1$ and using $E(u^*|Y,\Delta=1) = 0$ as required by (3.24),

$$-E[\mathbf{K}u^*|Y,\Delta=1] = E[\pi(Y,1,V)|Y,\Delta=1]f_{u^*}(Y)$$
$$+ E[\pi(Y,1,V)(Z - E(Z|Y,\Delta=1))|Y,\Delta=1].$$

Solving this for $f_{u^*}(Y)$ yields

(3.32)
$$f_{u^*}(Y) = -\frac{E[\mathbf{K}u^*|Y,\Delta=1]}{E[\pi(Y,1,V)|Y,\Delta=1]}$$
$$-\frac{E[\pi(Y,1,V)(Z - E(Z|Y,\Delta=1))|Y,\Delta=1]}{E[\pi(Y,1,V)|Y,\Delta=1]}.$$

Note that

(3.33)
$$\frac{E[\pi(Y,1,V)(Z - E(Z|Y,\Delta=1))|Y,\Delta=1]}{E[\pi(Y,1,V)|Y,\Delta=1]}$$
$$= \frac{E[\pi(Y,1,V)Z|Y,\Delta=1]}{E[\pi(Y,1,V)|Y,\Delta=1]} - E[Z|Y,\Delta=1]$$
$$= E[Z|Y,\Delta=1, R=1] - E[Z|Y,\Delta=1].$$

Combining (3.32) and (3.33) with (3.31), we obtain (3.21).  □

The reason that we solve for $u^*$ instead of $s^*$ is that there is no unique solution if we solve for $s^*$: if $s^*(Y,Z)$ satisfies (3.23) with $u^*$ replaced by $s^* - E[s^*|Y,\Delta=1]$, then $s^*(Y,Z) + f(Y)$ will satisfy (3.23) for any function $f(Y)$. From the form of $k \in \mathcal{K}$ we know that $k$ is determined by $u(Y,Z) \equiv s(Y,Z) - E[s(Y,Z)|Y,\Delta=1]$. So we only need to solve for $u^*(Y,Z) \equiv s^*(Y,Z) - E[s^*(Y,Z)|Y,\Delta=1]$ and then compute $l^*_\theta \equiv k^*$.

The parts of (3.18) and (3.21) involving $\mathbf{T}$ and $\mathbf{K}$, respectively, can be viewed as integral operators on the unknown function $u^*$. But these equations are not standard Fredholm integral equations of the second kind [see, e.g., Kress (1999) and Rudin (1973), Chapter 4]. Note that the terms of $\mathbf{K}$ involving conditional expectation given $Y' > Y$ correspond to noncompact operators [see, e.g., Rudin (1973), Problem 17, page 107] in general.

REMARK.  The equation corresponding to our (3.18) given by Robins, Rotnitzky and Zhao [(1994) Section 8.3, page 862, column 1] is incorrect. According to personal communications with Robins, their incorrect result is due to algebraic errors.

3.3. *Alternative models.* The survival model we have discussed so far assumes that $Y$ is always observable and $V$ is also a vector of covariates of interest. We now consider some alternative scenarios for the models involved in the two-phase designs.

*When $Z = (X, V)$ but the Cox model only involves $X$.* When $Z = (X, V)$, but the Cox model for the conditional distribution of $T$ given $Z$ includes only $X$, then $\theta'z$ is replaced by $\theta'x$ in the model (2.1). After going through the same procedure as deriving equation (3.18), we obtain

$$
\begin{aligned}
u^*(Y, Z) - \Big\{ &\mathbf{K}u^*(Y, Z) \\
&- \frac{\pi(Y, 1, V)}{E[\pi(Y, 1, V)|Y, \Delta = 1]} E[\mathbf{K}u^*(Y, Z)|Y, \Delta = 1] \Big\} \\
= \pi(Y, 1, V)\{X &- E[X|Y, R\Delta = 1]\},
\end{aligned}
\tag{3.34}
$$

where the operator $\mathbf{K}$ is as defined in (3.22).

*When $Y$ is not observed at phase 1.* If we are not able to observe $Y$ in the first phase, we will observe the data

$$
\begin{aligned}
(Y, \Delta, X, V)^R (\Delta, V)^{1-R} \\
\equiv (RY, \Delta, RX, V, R) \equiv \begin{cases} (Y, \Delta, X, V), & R = 1 \\ (\Delta, V), & R = 0, \end{cases}
\end{aligned}
$$

and the efficient score will have the form

$$
k^* = \frac{R}{\pi} \mathbf{D} u^*(Y, X, V) - \frac{R - \pi}{\pi} E[\mathbf{D} u^*(Y, X, V)|\Delta, V].
$$

By the same method used to derive (3.18), we find, for estimation of the coefficient of $Z$, that $u^*$ satisfies

$$
u^*(Y, Z) - \Big\{ \mathbf{K}u^*(Y, Z) - \frac{\pi(1, V)}{E[\pi(1, V)|Y, \Delta = 1]} E[\mathbf{K}u^*(Y, Z)|Y, \Delta = 1] \Big\}
= \pi(1, V)\{Z - E[Z|Y, R\Delta = 1]\},
\tag{3.35}
$$

where

$$
\begin{aligned}
\mathbf{K}u^*(Y, Z) = {}&-E[\mathbf{D}u^*(Y', \Delta, Z)|Y' > Y, Z] \\
&+ \pi(1, V) E\Big(\frac{\mathbf{D}u^*(Y', \Delta, Z)}{\pi(\Delta, V)}\Big|Y' > Y, Z\Big) \\
&+ (1 - \pi(1, V)) E[\mathbf{D}u^*|\Delta, V]|_{\Delta = 1} \\
&- \pi(1, V) E\Big(\frac{1 - \pi(\Delta, V)}{\pi(\Delta, V)} E[\mathbf{D}u^*(Y', \Delta, Z)|\Delta, V]\Big|Y' > Y, Z\Big).
\end{aligned}
\tag{3.36}
$$

*When $Y$ is not observed at phase 1, $Z = (X, V)$, and the Cox model only involves $X$.* When the Cox model involves only $\theta'X$ and just $(\Delta, V)$ is observed in the first phase, $u^*$ must satisfy

$$
\begin{aligned}
(3.37) \quad & u^*(Y, Z) - \left\{ \mathbf{K}u^*(Y, Z) - \frac{\pi(1, V)}{E[\pi(1, V)|Y, \Delta = 1]} E[\mathbf{K}u^*(Y, Z)|Y, \Delta = 1] \right\} \\
& = \pi(1, V)\{X - E[X|Y, R\Delta = 1]\},
\end{aligned}
$$

with $\mathbf{K}$ defined in (3.36).

*When $Z = (X, V)$, $V = (V_1, V_2)$ and the Cox model only involves $(X, V_1)$.* Suppose that $V = (V_1, V_2)$ and that $T$ and $V_2$ are conditionally independent given $(X, V_1)$, that is, the Cox model for the conditional distribution of $T$ involves only $(X, V_1)$. This is a generalization of the previous models. In order to avoid repeating calculation of the efficient score function, we can reduce the general model to the model either in Section 3.1 or in Section 3.3 and apply those existing results. Let $\tilde{X} = (X, V_1)$ and assume that $\tilde{X}$ is missing at random. Thus $\tilde{X}$ and $V$ would, respectively, play the role of $X$ and $V$ in our Sections 3.1 and 3.2. Although $V_1$ can be missing in $\tilde{X}$, it is fully recovered from $V$. Then we can directly use the results of alternative models in Sections 3.1 and 3.3 for the general model.

**4. Examples of information bound calculations.** The case-cohort design, studied by Prentice (1986) and Self and Prentice (1988), and the exposure stratified case-cohort design, studied by Borgan, Langholz, Samuelsen, Goldstein and Pogoda (2000), are two special cases in the class of two-phase designs. In the case-cohort design the complete information is essentially observed for all the failures and a simple random subsample of the nonfailures. The exposure stratified case-cohort design is a modification of the classical case-cohort design in which complete covariate data is observed for all failures and for a stratified random subsample of the nonfailures. The stratification is based upon a correlate (or surrogate variable, available for everyone) of the true exposure (or prognostic factor) of interest. In this section we treat the simplified i.i.d. versions of these sampling designs.

Pseudo-likelihood type (inefficient) estimators have been proposed by Prentice (1986) for case-cohort designs, and by Borgan, Langholz, Samuelsen, Goldstein and Pogoda (2000) for exposure stratified case-cohort designs. For discussions of efficient estimators for these designs we refer to Nan (2001). But information bound calculations can tell us how much information we could potentially gain from fully efficient estimators and which design methods use the observed data more efficiently, if efficient estimators were available. Here we give two examples in which the information bound calculations can be carried out analytically. Although these two hypothetical examples are rather special cases involving evaluation of the information for a simple parametric subfamily, the calculations in this section may tell us the fundamental properties of the designs and potential estimators. The results may give some guidance for designing and analyzing real studies.

4.1. *A case-cohort study.* We assume that the true distribution has exponentially distributed failure times and a single binary covariate $Z$ taking values 0 and 1. Let $h(z) = P(Z = z)$ be the probability that a subject has covariate value $z \in \{0, 1\}$; thus $h(0) + h(1) = 1$. The censoring time is distributed with point mass 1 at $t = 1$, which means that all subjects in the cohort are followed from time zero to either failure or to the end of the study at $t = 1$. This is discussed as an example by Self and Prentice (1988). The density of the complete data can be written as

$$(4.1) \quad q_{Y,\Delta,Z}(y,\delta,z) = \begin{cases} w_1(y|z)h(z) = \lambda e^{\theta z - \lambda y e^{\theta z}} h(z), & \delta = 1,\ 0 \leq y \leq 1, \\ w_2(y|z)h(z) = e^{-\lambda e^{\theta z}} h(z), & \delta = 0,\ y = 1. \end{cases}$$

In the i.i.d. version of the case-cohort study, a simple random subsample is taken from the nonfailures with sampling (inclusion) probability $\pi_0$. The values of the covariate $Z$ are measured for all the failures and only the sampled nonfailures. Hence it is a special case of the two-phase design discussed in the previous section with $\pi(Y, \Delta) \equiv P(R = 1 | Y, \Delta) = \pi(\Delta)$ only, and where $\pi(1) = 1$, $\pi(0) = \pi_0 \in (0, 1)$. Note that we do not have a surrogate covariate $V$ in this example. In a classical case-cohort design, $Y$ may not be observed if the subject is not a failure and not in the subcohort. But for this special example $\Delta = 0$ implies $Y = 1$. So for information bound calculations it does not matter whether we treat $Y$ as known or not. Detailed calculation is omitted here and can be found in Nan (2001).

Figure 1 displays the ratios of asymptotic variance of the Self and Prentice (1988) pseudo-likelihood estimator (SP Variance) to the information lower bound for $\theta$ as a function of the sampling fraction for nonfailures in the i.i.d. case-cohort model shown above. Figure 1 shows that when the disease is rare, that is, the baseline failure probability is very low, the pseudo-likelihood estimator is close to fully efficient. As the failure probability increases, the pseudo-likelihood estimator loses more efficiency, especially when the subcohort fraction is small. Hence development of more efficient estimators may be worthwhile for case-cohort designs where increasing the subcohort fraction is costly and the failure probability is moderate.

Figure 2 displays the ratio of the information lower bound for estimation of $\theta$ based on the "observed data" ($1/I_\theta^*$ where $I_\theta^*$ is the information for $\theta$) and the asymptotic variance of the partial-likelihood estimator for $\theta$ based on "complete" (or "full") data. This ratio is shown as a function of the sampling fraction for the nonfailures under different baseline failure probabilities when

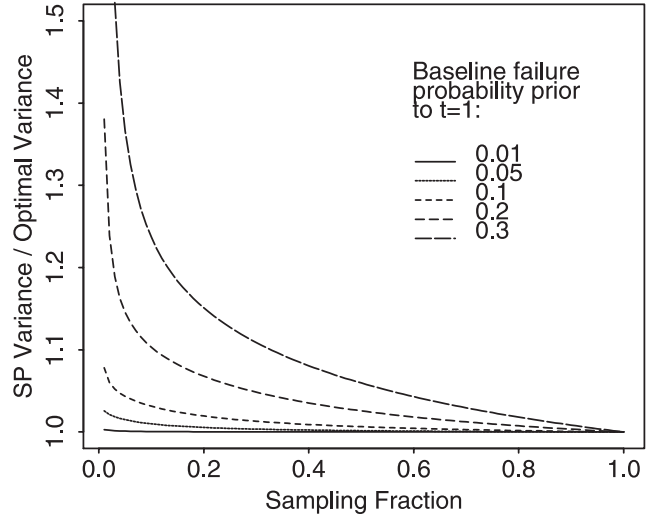

FIG. 1. *Ratio of the variance of the Self and Prentice pseudo-likelihood estimator* (SP) *to the Optimal Variance as a function of the sampling fraction for nonfailures under different baseline failure probabilities, $P(T \leq 1|Z=0)$. Here $\theta = \ln(2)$.*

$e^\theta = 2$. Figure 2 shows that the case-cohort design loses more information (supposing that an efficient estimator is available), relative to complete data, as the failure probability increases and as the subcohort fraction decreases.

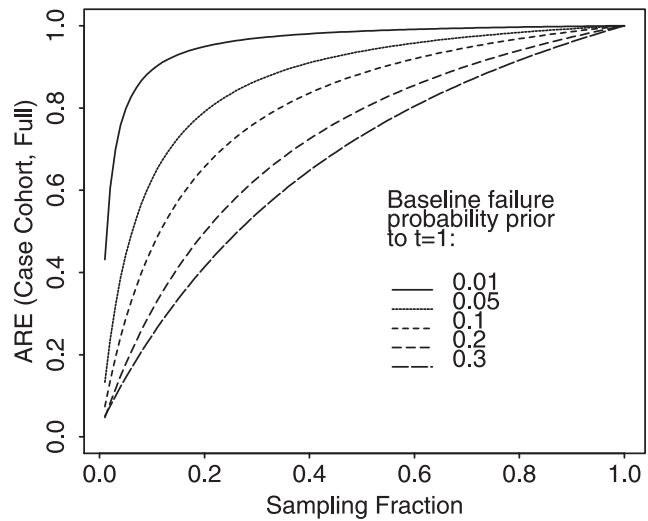

FIG. 2. *Asymptotic relative efficiency of optimally efficient estimators for the i.i.d. case-cohort design as a function of the sampling fraction for nonfailures under different baseline failure probabilities, $P(T \leq 1|Z=0)$. Here $\theta = \ln(2)$.*

Actually, when the baseline failure probability is above 0.5, the curves in Figure 2 move toward the upper left again as the failure probability increases, but there is less interest in these high failure probability cases in practice. From Figure 2 we can see that a great deal of precision may be lost by using a case-cohort design as opposed to data collection on the full cohort, even when a fully efficient estimator is used for the case-cohort study. With this knowledge investigators can weigh the trade-off between precision and study cost. Further work is needed to explore presumably more efficient designs: for example, an alternative design might be an "exposure stratified case-cohort design" as in our second example.

Perhaps the more interesting phenomena appear in Figure 3 and Figure 4. In Figure 3 we look at the asymptotic relative efficiency of the pseudo-likelihood estimator as a function of $\theta$. Figure 4 shows the relative efficiency of the optimal variances for the i.i.d. case-cohort design versus the full data design as a function of $\theta$. When $\theta$ is near zero Figure 3 shows that the pseudo-likelihood estimator does not lose much efficiency compared to the optimal estimator for the case-cohort design. However, Figure 4 shows that the case-cohort design (with an optimal estimator) loses considerable information compared to the full data design. The minimum ARE (as a function of $\theta$) depends on the baseline failure probability; the minimum increases and it moves away from $\theta = 0$ as the baseline failure probability decreases. When $\theta$ is away from zero, that is, the effect of the covariate $Z$ is large, the pseudo-likelihood estimator loses significant efficiency, especially when $\theta$ is positive and the baseline failure probability is high. However, away from zero the design itself starts to gain information and is very close to the full data design when the absolute value of $\theta$ is large. The conclusion is that if we expect intermediate to large covariate effects, it may be very worthwhile to find efficient estimators for $\theta$. Certainly, developing more efficient designs is also valuable, as can be seen from Figures 2 and 4.

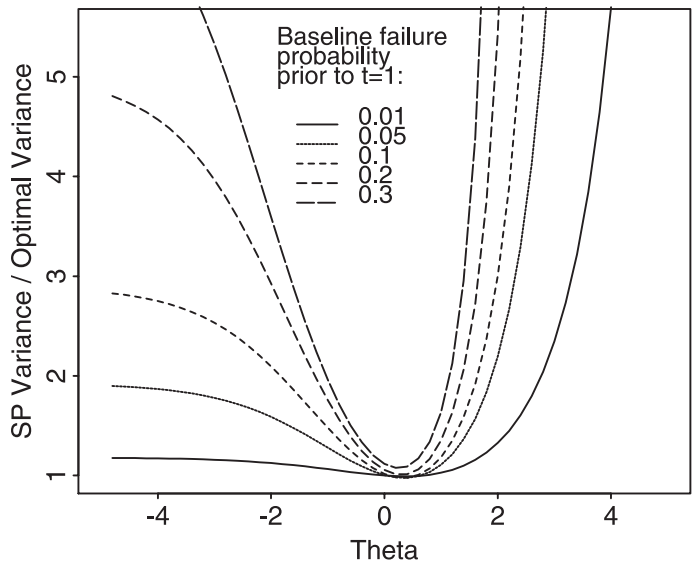

FIG. 3. *Ratio of the variance for the Self and Prentice pseudo-likelihood estimator* (SP) *to the Optimal Variance as a function of* $\theta$ (*log of relative risk*) *for different baseline failure probabilities* $P(T \leq 1 | Z = 0)$ *in the i.i.d. case-cohort design. Here the sampling fraction for nonfailures is* 0.1.

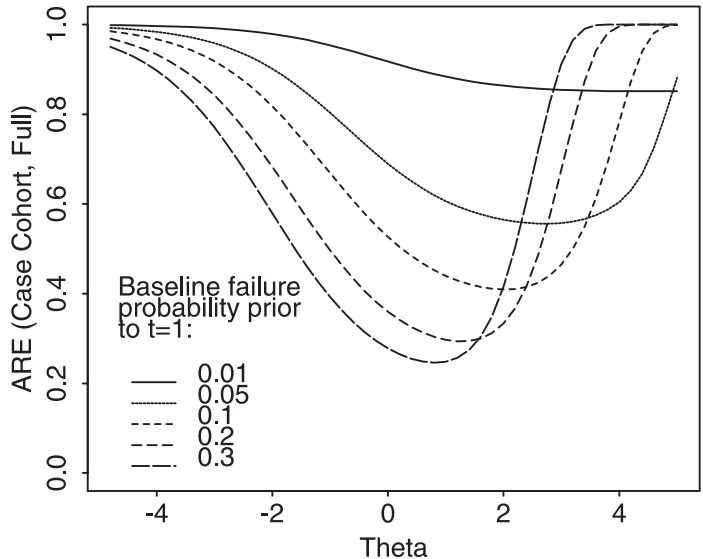

FIG. 4. *Asymptotic relative efficiency of the optimally efficient estimator in the case-cohort design relative to the optimally efficient estimator in the full cohort design as a function of $\theta$ (log of relative risk) for different baseline failure probabilities $P(T \leq 1|Z=0)$. Here the sampling fraction for nonfailures is* 0.1.

4.2. *An exposure stratified case-cohort study.* Assume that $X$ is the variable of interest and that $V$ is a surrogate variable for $X$, or measurement of $X$ with error, and $V$ is conditionally independent of $T$ given $X$. We suppose that $V$ can be observed for everyone in the entire cohort, but $X$ is only observed for subjects in the subcohort and failures. Then the model for this type of data is the first alternative model discussed in the previous section. The i.i.d. version of the exposure stratified case-cohort design studied by Borgan, Langholz, Samuelsen, Goldstein and Pogoda (2000) is a special case of this model. Here we discuss an example with a binary covariate $X \in \{0,1\}$ and a binary surrogate variable $V \in \{0,1\}$. The distribution for $X$ and $V$ has the form of a $2 \times 2$ table. Let

$$P(V=1|X=1) = 1-\alpha, \qquad P(V=0|X=0) = 1-\beta.$$

If we consider $V = 1$ as a "positive" test for $X$, then $1 - \alpha$ is the sensitivity and $1 - \beta$ is the specificity of the test. We assume exponentially distributed failure times. All subjects in the cohort are followed from time zero to either failure or to the end of the study at time $t = 1$. For our calculations the exponential failure rate parameter $\lambda$ will be set to achieve a specified baseline failure probability as in Section 4.1. Let the joint mass function of $(X, V)$ be $h(x, v)$. Thus we have the joint density for the underlying complete data,

$$(4.2) \quad q_{Y,\Delta,X,V}(y, \delta, x, v) = \begin{cases} \lambda e^{\theta x - \lambda y e^{\theta x}} h(x, v), & \delta = 1,\ 0 \leq y \leq 1, \\ e^{-\lambda e^{\theta x}} h(x, v), & \delta = 0,\ y = 1. \end{cases}$$

By the same argument as in the previous example, we may assume that $Y$ always is observed. The cohort is then categorized into three strata: $\{\Delta = 1\}$, $\{\Delta = 0, V = 0\}$ and $\{\Delta = 0, V = 1\}$. We observe complete information for all the subjects in the first stratum, and of $\pi_0$, $\pi_1$ fractions (constants) of the subjects in the second and third strata, respectively. We only observe $(Y, \Delta, V)$ for other subjects. In probability language we have $P(R = 1|\Delta = 1, Y, V) \equiv \pi(Y, 1, V) = 1$, $P(R = 1|\Delta = 0, Y, V) \equiv \pi(Y, 0, V) \equiv \pi_0$ if $V = 0$ and $\pi_1$ if $V = 1$. Again, we omit the detailed calculations here and refer to Nan (2001).

We calculate $I_\theta^*$ for different $\alpha$, $\beta$, $P(X = 0)$, $\pi_0$, $\pi_1$, $\theta$ and $\lambda$ by using numerical integration. When $\alpha = \beta = 0.5$, the exposure stratified case-cohort

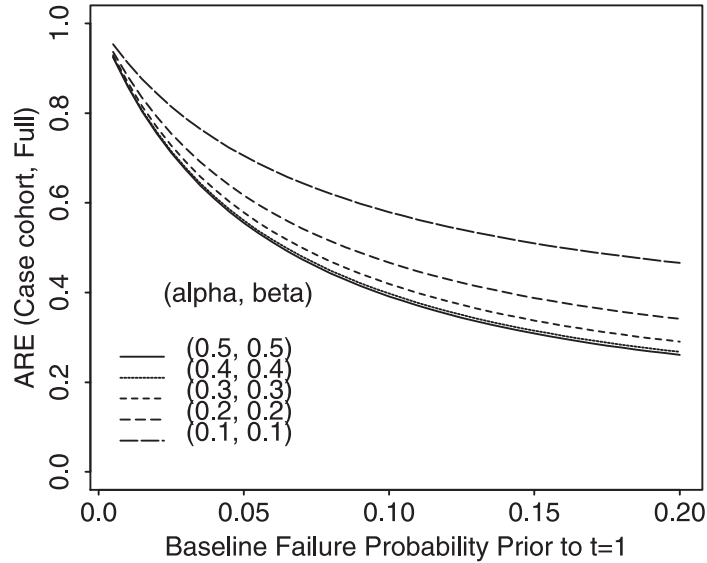

FIG. 5. *Asymptotic efficiency of the optimally efficient estimator for the case-cohort design with a surrogate variable, relative to that for the full cohort design, as a function of the baseline failure probability, $P(T \leq 1|X = 0)$. Here $\theta = \ln(2)$, $P(X = 0) = 0.9$, $\pi_0 = \pi_1 = 0.1$ (i.e., stage 2 sampling is not stratified).*

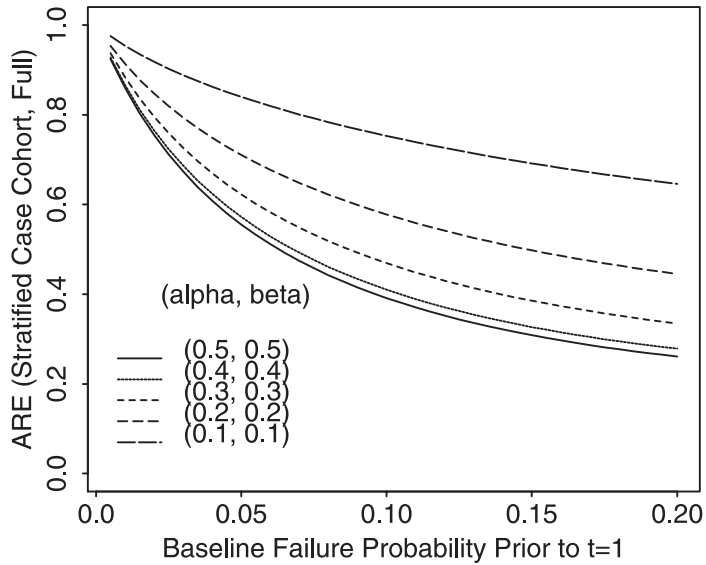

Fig. 6. *Asymptotic efficiency of the optimally efficient estimator for the i.i.d. stratified case-cohort design with a surrogate variable, relative to that for the full cohort design, as a function of the baseline failure probability $P(T \leq 1 | X = 0)$. Here, $\theta = \ln(2)$, $P(X = 0) = 0.9$, $\pi_0 \neq \pi_1$ (i.e., with stratified sampling at stage 2).*

design is equivalent to the classical case-cohort design (previous example) since $V$ is not correlated with $X$ under this condition. Figures 5 and 6 show the comparisons of the asymptotic relative efficiency (ARE) of fully efficient estimators (if they exist) for the exposure stratified (at $\pi_0 = \pi_1$ and $\pi_0 \neq \pi_1$) and classical case-cohort designs at $e^\theta = 2$ and $P(X = 0) = 0.9$. When $\theta = 0$, the corresponding figures (not shown) have similar patterns but slightly different magnitude. In Figure 5 the sampling probabilities in the two strata are equal, that is, $\pi_0 = \pi_1 = 0.1$. In Figure 6 $\pi_0$ and $\pi_1$ are different, such that the expected numbers of sampled subjects in strata $\{\Delta = 0, V = 0\}$ and $\{\Delta = 0, V = 1\}$ are the same (or approximately the same), and the fraction of sampled subjects from the two strata all together is 0.1 (or approximately 0.1). We can see from Figure 5 that the efficiency increases as the sensitivity $(1 - \alpha)$ and specificity $(1 - \beta)$ increase as a result of incorporating the surrogate variable into the model, even without stratification. Note that an efficient estimator will incorporate $V$ for subjects outside the subsample, providing information for the estimation of $\theta$ via the correlation between $X$ and $V$. When we do stratified sampling $(\pi_0 \neq \pi_1)$, Figure 6 shows that the efficiency gains are even greater. So both incorporating surrogate information and stratified sampling will increase the efficiency. Note that the information bound calculation illustrated here is for a measurement er-

TABLE 1
*Comparisons of asymptotic relative efficiency (ARE relative to a full cohort study): the approximate AREs of pseudo-likelihood estimators and the AREs of information bounds for a stratified case-cohort design, which has only the binary covariate of interest, $X$, and a binary covariate $V$, which is a surrogate for $X$. The specificity of $V = 1$ as a test for $X = 1$ is $1 - \beta$, and the sensitivity is $1 - \alpha$. The subcohort size equals the expected number of cases ($PL =$ Pseudo-Likelihood. This part is taken from Table 1 of Borgan, Langholz, Samuelsen, Goldstein anf Pogoda (2000). $IB =$ Information Bound.)*

(a) $P(X = 1) = 0.05$

|  | $ARE(PL), \%$ | | | $ARE(IB), \%$ | | | $\frac{ARE(PL)}{ARE(IB)}, \%$ | | |
|---|---|---|---|---|---|---|---|---|---|
|  | $1 - \beta$ | | | $1 - \beta$ | | | $1 - \beta$ | | |
| $1 - \alpha$ | 0.50 | 0.70 | 0.90 | 0.50 | 0.70 | 0.90 | 0.50 | 0.70 | 0.90 |
| 0.50 | 35.5 | 36.5 | 40.8 | 36.0 | 37.8 | 45.9 | 98.6 | 96.6 | 88.9 |
| 0.70 | 36.5 | 39.6 | 47.3 | 37.7 | 43.0 | 55.9 | 96.8 | 92.1 | 84.6 |
| 0.90 | 40.8 | 47.3 | 60.5 | 43.9 | 53.0 | 70.3 | 92.9 | 89.2 | 86.1 |

(b) $P(X = 1) = 0.50$

|  | $1 - \beta$ | | | $1 - \beta$ | | | $1 - \beta$ | | |
|---|---|---|---|---|---|---|---|---|---|
| $1 - \alpha$ | 0.50 | 0.70 | 0.90 | 0.50 | 0.70 | 0.90 | 0.50 | 0.70 | 0.90 |
| 0.50 | 52.9 | 54.0 | 58.4 | 53.5 | 55.5 | 63.2 | 98.9 | 97.3 | 92.4 |
| 0.70 | 54.0 | 57.3 | 64.7 | 55.5 | 61.1 | 72.0 | 97.3 | 93.8 | 89.9 |
| 0.90 | 58.4 | 64.7 | 75.8 | 62.6 | 71.2 | 83.6 | 93.3 | 90.9 | 90.7 |

ror problem without making any assumptions on the structure of the joint distribution of $(X, V)$.

Borgan, Langholz, Samuelsen, Goldstein and Pogoda (2000) study inverse probability weighted estimators in this model. In order to show how much efficiency the inverse probability estimator loses, we calculate the ARE of their estimator relative to a fully efficient estimator with asymptotic variance given by our information bound $1/I_\theta^*$ for the setting of their Example 1. We choose their optimal sampling fractions and a very small failure rate, $\lambda = 0.01$, which we believe is small enough to be able to make valid comparisons to their results. Note that the results in Table 1 are calculated under the condition that the subcohort fraction equals the expected number of cases, providing approximately one "control" per case, a frequently used design.

**5. Conclusions and further problems.** We have established new information bounds for the Cox model with missing data. Along the way we have developed a new decomposition of $L_2^0(Q)$, characterized the structure of the orthogonal complement of the nuisance parameter tangent space $\dot{\mathcal{Q}}_\eta^\perp$, and shown how to project onto the space $\dot{\mathcal{Q}}_\eta^\perp$ using conditional versions $R_1$ and $R_2$ of the mean residual life operator $R$ introduced by Efron and Johnstone (1990) and Ritov and Wellner (1988). The new bounds can be used to examine the loss of efficiency of pseudo-likelihood estimators for a given design and the amount of information loss due to a given design relative to complete data collection or to an alternative two phase design. While it has been known for some time that pseudo-likelihood estimators are not semiparametrically efficient, our explicit calculations quantify the loss of efficiency and also show that two-phase designs with stratified subsampling can partially recover the information that is lost due to missing data.

*Further problems.*

1. *Construction of efficient estimators when covariates are discrete.* Efficient estimators can be constructed explicitly using one-step methods when the covariates are discrete. For a preliminary study of such estimators, see Nan (2001).
2. *Construction of efficient estimators in general.* This will depend crucially on understanding the properties of the integral equation defining the efficient score and influence function. A major difficulty in constructing efficient estimators is the fact that the conditional cumulative hazard function $\Lambda_G(y|z)$ enters into the key equation (3.18) which determines $u^*$ and hence the efficient score function $l_\theta^*$. This function is typically completely unknown and is a function of $d+1$ variables which must be estimated nonparametrically. This is, of course, a difficult task for even moderately large $d$. However, our goal is not to estimate $\Lambda_G$ well, but instead to estimate $\theta$ well, and it is not yet clear how crucial the difficulty in estimating $\Lambda_G$ will be for construction of (nearly) efficient estimates of $\theta$. We remain optimistic about this at least for moderate values of $d$, and regard this as an important question for future work.
3. *How can we "optimize" the sampling design for a particular study?* If we focus on the variance of the estimator of a particular regression coefficient (e.g., the coefficient corresponding to a binary treatment-control covariate), then it would be very interesting to know how to allocate the sampling effort in the second phase to minimize the (asymptotic) variance. Our results provide the tools to graphically address this extremely important question.
4. *Are there better compromise estimators based on pseudo-likelihood?* Here the approaches of Chatterjee, Chen and Breslow (2003) and Chatterjee (1999) may be useful.

**Acknowledgments.** We owe thanks to Norman Breslow, Nilanjan Chatterjee and the other participants in the *Missing Data Working Group* at the University of Washington during the period 1995–1998, for many useful discussions about missing data and the subject of this paper. We also owe thanks to J. Robins for helpful discussions and to a referee for suggesting the current form given in Proposition 3.3 for the equation characterizing $u^*$.

B. Nan  
Department of Biostatistics  
University of Michigan  
1420 Washington Heights  
Ann Arbor, Michigan 48109-2029  
USA  
e-mail: bnan@umich.edu

M. Emond  
Department of Biostatistics  
University of Washington  
P.O. Box 357232  
Seattle, Washington 98195-7232  
USA  
e-mail: emond@u.washington.edu

J. A. Wellner  
Department of Statistics  
University of Washington  
P.O. Box 354322  
Seattle, Washington 98195-4322  
USA  
e-mail: jaw@stat.washington.edu